\documentclass[12pt]{amsart}
\usepackage[utf8]{inputenc}

\usepackage{subfiles}
\usepackage{hyperref}
\usepackage{cleveref}
\usepackage{amsmath}
\usepackage{amssymb}
\usepackage{amsthm}
\usepackage{enumitem}
\usepackage{mathtools}
\usepackage{thmtools}
\usepackage{thm-restate}
\usepackage[margin=1in]{geometry}
\usepackage{cleveref}
\usepackage{ytableau}
\usepackage{comment}
\usepackage{soul}
\usepackage{tikz}
\usetikzlibrary{tikzmark,decorations.pathreplacing}
\usetikzlibrary{arrows,matrix}
\usetikzlibrary{positioning}
\usetikzlibrary{decorations}
\usepackage{cite}
\usepackage{graphicx}
\usepackage{caption}
\usepackage{float}

\numberwithin{equation}{section}

\let\oldbibliography\thebibliography
\renewcommand{\thebibliography}[1]{
  \oldbibliography{#1}
  \setlength{\itemsep}{2pt}
}

\theoremstyle{definition}
\newtheorem{theorem}{Theorem}[section]
\newtheorem*{theorem*}{Theorem}

\newtheorem*{example*}{Example}
\newtheorem{lemma}[theorem]{Lemma}
\newtheorem*{lemma*}{Lemma}
\newtheorem{corollary}[theorem]{Corollary}
\newtheorem*{corollary*}{Corollary}

\newtheorem*{definition*}{Definition}
\newtheorem{proposition}[theorem]{Proposition}
\newtheorem*{proposition*}{Proposition}

\newtheorem*{remark*}{Remark}
\newtheorem{conjecture}[theorem]{Conjecture}

\usepackage{ mathrsfs }

\title{Stack-Sorting with Dotted-Pattern-Avoiding Stacks}

\author{Michael Yang}\address{\textsc{M. Yang},
Phillips Exeter Academy,
    Exeter, MA, 03833}
    \email{michaelyang800@gmail.com}
    
\author{Hansen Shieh}\address{\textsc{H. Shieh}, Westford Academy,
    Westford, MA, 01886} \email{cyborgpenguinto@gmail.com}
    
\author{Ashley Yu}\address{\textsc{A. Yu}, Concord Academy,
    Concord, MA, 01742} \email{ashley.y.ca99@gmail.com}
\begin{document}

\begin{abstract}
In this paper, we introduce the dotted pattern-avoiding map $s_{\dot{\tau}}$, which avoids the dotted pattern $\dot{\tau}$ instead of descents as West's stack-sorting map $s$ does. We also extend the pattern-avoiding machine, which is composed of a $\sigma$-avoiding map and West's stack-sorting map $s$, to the dotted pattern-avoiding machine. In this paper, we prove the analogs of the classical results on West's stack-sorting map for the length-$2$ dotted pattern-avoiding maps. We end with several conjectures. 
\end{abstract}

\maketitle

\section{Introduction} 
\label{intro}
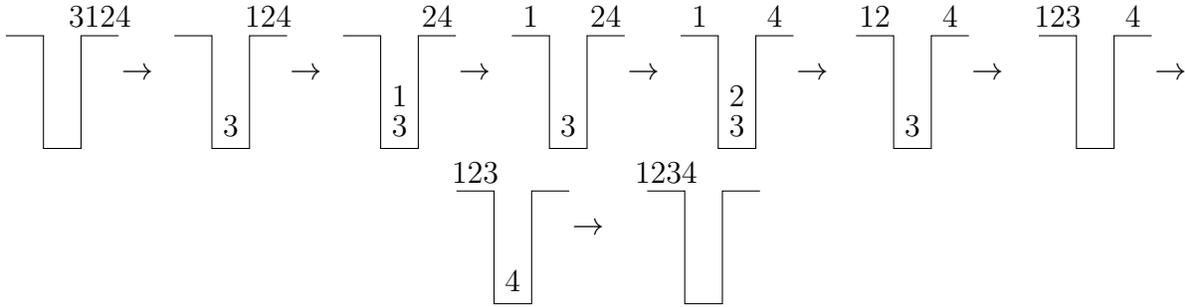
\begin{figure}
\begin{center}
\begin{tikzpicture}[scale=0.5]
    \draw (0,0) -- (1,0) -- (1,-3) -- (2,-3) -- (2,0) -- (3,0);
    \node[fill = none, draw = none] at (2.5,.5) {3124};
    \node[fill = none, draw = none] at (1.5,-2.4) {};
    \node[fill = none, draw = none] at (1.5,-1.6) {};
    \node[fill = none, draw = none] at (1.5,-.8) {};
    \node[fill = none, draw = none] at (3.5,-1) {$\rightarrow$};
\end{tikzpicture}
\begin{tikzpicture}[scale=0.5]
    \draw (0,0) -- (1,0) -- (1,-3) -- (2,-3) -- (2,0) -- (3,0);
    \node[fill = none, draw = none] at (2.5,.5) {124};
    \node[fill = none, draw = none] at (1.5,-2.4) {3};
    \node[fill = none, draw = none] at (1.5,-1.6) {};
    \node[fill = none, draw = none] at (1.5,-.8) {};
    \node[fill = none, draw = none] at (3.5,-1) {$\rightarrow$};
\end{tikzpicture}
\begin{tikzpicture}[scale=0.5]
    \draw (0,0) -- (1,0) -- (1,-3) -- (2,-3) -- (2,0) -- (3,0);
    \node[fill = none, draw = none] at (2.5,.5) {24};
    \node[fill = none, draw = none] at (1.5,-2.4) {3};
    \node[fill = none, draw = none] at (1.5,-1.6) {1};
    \node[fill = none, draw = none] at (1.5,-.8) {};
    \node[fill = none, draw = none] at (3.5,-1) {$\rightarrow$};
\end{tikzpicture}
\begin{tikzpicture}[scale=0.5]
    \draw (0,0) -- (1,0) -- (1,-3) -- (2,-3) -- (2,0) -- (3,0);
    \node[fill = none, draw = none] at (2.5,.5) {24};
    \node[fill = none, draw = none] at (.5,.5) {1};
    \node[fill = none, draw = none] at (1.5,-2.4) {3};
    \node[fill = none, draw = none] at (1.5,-1.6) {};
    \node[fill = none, draw = none] at (1.5,-.8) {};
    \node[fill = none, draw = none] at (3.5,-1) {$\rightarrow$};
\end{tikzpicture}
\begin{tikzpicture}[scale=0.5]
    \draw (0,0) -- (1,0) -- (1,-3) -- (2,-3) -- (2,0) -- (3,0);
    \node[fill = none, draw = none] at (2.5,.5) {4};
    \node[fill = none, draw = none] at (.5,.5) {1};
    \node[fill = none, draw = none] at (1.5,-2.4) {3};
    \node[fill = none, draw = none] at (1.5,-1.6) {2};
    \node[fill = none, draw = none] at (1.5,-.8) {};
    \node[fill = none, draw = none] at (3.5,-1) {$\rightarrow$};
\end{tikzpicture}
\begin{tikzpicture}[scale=0.5]
    \draw (0,0) -- (1,0) -- (1,-3) -- (2,-3) -- (2,0) -- (3,0);
    \node[fill = none, draw = none] at (2.5,.5) {4};
    \node[fill = none, draw = none] at (.5,.5) {12};
    \node[fill = none, draw = none] at (1.5,-2.4) {3};
    \node[fill = none, draw = none] at (1.5,-1.6) {};
    \node[fill = none, draw = none] at (1.5,-.8) {};
    \node[fill = none, draw = none] at (3.5,-1) {$\rightarrow$};
\end{tikzpicture}
\begin{tikzpicture}[scale=0.5]
    \draw (0,0) -- (1,0) -- (1,-3) -- (2,-3) -- (2,0) -- (3,0);
    \node[fill = none, draw = none] at (2.5,.5) {4};
    \node[fill = none, draw = none] at (.5,.5) {123};
    \node[fill = none, draw = none] at (1.5,-2.4) {};
    \node[fill = none, draw = none] at (1.5,-1.6) {};
    \node[fill = none, draw = none] at (1.5,-.8) {};
    \node[fill = none, draw = none] at (3.5,-1) {$\rightarrow$};
\end{tikzpicture}
\begin{tikzpicture}[scale=0.5]
    \draw (0,0) -- (1,0) -- (1,-3) -- (2,-3) -- (2,0) -- (3,0);
    \node[fill = none, draw = none] at (2.5,.5) {};
    \node[fill = none, draw = none] at (.5,.5) {123};
    \node[fill = none, draw = none] at (1.5,-2.4) {4};
    \node[fill = none, draw = none] at (1.5,-1.6) {};
    \node[fill = none, draw = none] at (1.5,-.8) {};
    \node[fill = none, draw = none] at (3.5,-1) {$\rightarrow$};
\end{tikzpicture}
\begin{tikzpicture}[scale=0.5]
    \draw (0,0) -- (1,0) -- (1,-3) -- (2,-3) -- (2,0) -- (3,0);
    \node[fill = none, draw = none] at (2.5,.5) {};
    \node[fill = none, draw = none] at (.5,.5) {1234};
    \node[fill = none, draw = none] at (1.5,-2.4) {};
    \node[fill = none, draw = none] at (1.5,-1.6) {};
    \node[fill = none, draw = none] at (1.5,-.8) {};
\end{tikzpicture}
\end{center}
\caption{West's stack-sorting map $s$ on the permutation $\pi=3124$}
\label{fig1}
\end{figure}

In 1968, Knuth \cite{knuth} first introduced the stack-sorting machine that, at each step, either \textit{pushes} the leftmost element of the input permutation to the stack or \textit{pops} the top element of the stack onto the output. In 1990, West \cite{west} introduced a deterministic stack-sorting map $s$ that drew inspiration from Knuth's non-deterministic stack-sorting machine. West's stack-sorting map $s$ processes the input permutation through a stack in a right greedy manner such that elements of the stack always increase from top to bottom (see for example, \Cref{fig1}). 

As hinted by the name, West's stack-sorting map eventually sorts all permutations after repeated iterations. A natural question, then, is how many permutations of length $n$ get sorted after $t$ iterations through West's stack-sorting map. 

As a first step to answering the question, West defined a permutation $\pi$ as \textit{$t-$stack-sortable} if $s^{t}(\pi)$ is the identity permutation. Subsequently, Knuth demonstrated that a permutation is 1-stack-sortable if and only if it does not contain a subsequence that is order-isomorphic to 231. Knuth moreover proved that the number of 1-stack-sortable permutations in $S_n$ is the $n^{\text{th}}$ Catalan number $C_n = \tfrac{1}{n + 1}\tbinom{2n}{n}$. In 1990, West \cite{west} proposed a conjecture that the number of 2-stack-sortable permutations in $S_n$ is $\frac{2}{(n+1)(2n+1)}\binom{3n}{n}$, which was subsequently proved by Zeilberger \cite{zeil}. 

In 2021, Defant \cite{defant} defined the dual of a $t-$stack-sortable permutation as a $t$-sorted permutation for $t$ that is close to $n$. A permutation is \textit{$t-$sorted} if it is in the image of $s^t(S_n)$. In a related direction, Zhang \cite{zhang} defined a \textit{minimally sorted} permutation of $S_n$ to be a permutation that requires the maximum number of sorts through the stack-sorting map to be sorted. 

Another chain of research extended West's stack-sorting map. In 2020, Cerbai, Claesson, and Ferrari \cite{stacksort} generalized $s$ to $s_{\sigma}$, in which the stack must avoid the pattern $\sigma$ rather than descents. In addition, Cerbai, Claesson, and Ferrari further extended pattern-avoiding maps to $\sigma$-machines, which are composed of $s_{\sigma}$ followed by $s$. Berlow \cite{berlow} extended $\sigma$-machines to maps avoiding two patterns and defined the $(\sigma,\tau)$-machine as the composition of $s_{\{\sigma,\tau\}}$ followed by $s$. More recently, Defant and Zheng \cite{colinandkai} extended the stack-sorting map to stacks that avoid consecutive patterns.
 
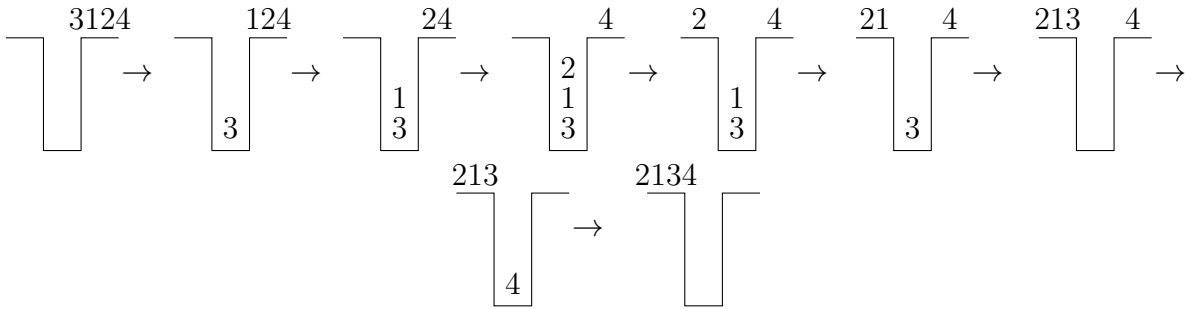
\begin{figure}[h]
\begin{center}
\begin{tikzpicture}[scale=0.5]
    \draw (0,0) -- (1,0) -- (1,-3) -- (2,-3) -- (2,0) -- (3,0);
    \node[fill = none, draw = none] at (2.5,.5) {3124};
    \node[fill = none, draw = none] at (1.5,-2.4) {};
    \node[fill = none, draw = none] at (1.5,-1.6) {};
    \node[fill = none, draw = none] at (1.5,-.8) {};
    \node[fill = none, draw = none] at (3.5,-1) {$\rightarrow$};
\end{tikzpicture}
\begin{tikzpicture}[scale=0.5]
    \draw (0,0) -- (1,0) -- (1,-3) -- (2,-3) -- (2,0) -- (3,0);
    \node[fill = none, draw = none] at (2.5,.5) {124};
    \node[fill = none, draw = none] at (1.5,-2.4) {3};
    \node[fill = none, draw = none] at (1.5,-1.6) {};
    \node[fill = none, draw = none] at (1.5,-.8) {};
    \node[fill = none, draw = none] at (3.5,-1) {$\rightarrow$};
\end{tikzpicture}
\begin{tikzpicture}[scale=0.5]
    \draw (0,0) -- (1,0) -- (1,-3) -- (2,-3) -- (2,0) -- (3,0);
    \node[fill = none, draw = none] at (2.5,.5) {24};
    \node[fill = none, draw = none] at (1.5,-2.4) {3};
    \node[fill = none, draw = none] at (1.5,-1.6) {1};
    \node[fill = none, draw = none] at (1.5,-.8) {};
    \node[fill = none, draw = none] at (3.5,-1) {$\rightarrow$};
\end{tikzpicture}
\begin{tikzpicture}[scale=0.5]
    \draw (0,0) -- (1,0) -- (1,-3) -- (2,-3) -- (2,0) -- (3,0);
    \node[fill = none, draw = none] at (2.5,.5) {4};
    \node[fill = none, draw = none] at (1.5,-2.4) {3};
    \node[fill = none, draw = none] at (1.5,-1.6) {1};
    \node[fill = none, draw = none] at (1.5,-.8) {2};
    \node[fill = none, draw = none] at (3.5,-1) {$\rightarrow$};
\end{tikzpicture}
\begin{tikzpicture}[scale=0.5]
    \draw (0,0) -- (1,0) -- (1,-3) -- (2,-3) -- (2,0) -- (3,0);
    \node[fill = none, draw = none] at (2.5,.5) {4};
    \node[fill = none, draw = none] at (.5,.5) {2};
    \node[fill = none, draw = none] at (1.5,-2.4) {3};
    \node[fill = none, draw = none] at (1.5,-1.6) {1};
    \node[fill = none, draw = none] at (1.5,-.8) {};
    \node[fill = none, draw = none] at (3.5,-1) {$\rightarrow$};
\end{tikzpicture}
\begin{tikzpicture}[scale=0.5]
    \draw (0,0) -- (1,0) -- (1,-3) -- (2,-3) -- (2,0) -- (3,0);
    \node[fill = none, draw = none] at (2.5,.5) {4};
    \node[fill = none, draw = none] at (.5,.5) {21};
    \node[fill = none, draw = none] at (1.5,-2.4) {3};
    \node[fill = none, draw = none] at (1.5,-1.6) {};
    \node[fill = none, draw = none] at (1.5,-.8) {};
    \node[fill = none, draw = none] at (3.5,-1) {$\rightarrow$};
\end{tikzpicture}
\begin{tikzpicture}[scale=0.5]
    \draw (0,0) -- (1,0) -- (1,-3) -- (2,-3) -- (2,0) -- (3,0);
    \node[fill = none, draw = none] at (2.5,.5) {4};
    \node[fill = none, draw = none] at (.5,.5) {213};
    \node[fill = none, draw = none] at (1.5,-2.4) {};
    \node[fill = none, draw = none] at (1.5,-1.6) {};
    \node[fill = none, draw = none] at (1.5,-.8) {};
    \node[fill = none, draw = none] at (3.5,-1) {$\rightarrow$};
\end{tikzpicture}
\begin{tikzpicture}[scale=0.5]
    \draw (0,0) -- (1,0) -- (1,-3) -- (2,-3) -- (2,0) -- (3,0);
    \node[fill = none, draw = none] at (2.5,.5) {};
    \node[fill = none, draw = none] at (.5,.5) {213};
    \node[fill = none, draw = none] at (1.5,-2.4) {4};
    \node[fill = none, draw = none] at (1.5,-1.6) {};
    \node[fill = none, draw = none] at (1.5,-.8) {};
    \node[fill = none, draw = none] at (3.5,-1) {$\rightarrow$};
\end{tikzpicture}
\begin{tikzpicture}[scale=0.5]
    \draw (0,0) -- (1,0) -- (1,-3) -- (2,-3) -- (2,0) -- (3,0);
    \node[fill = none, draw = none] at (2.5,.5) {};
    \node[fill = none, draw = none] at (.5,.5) {2134};
    \node[fill = none, draw = none] at (1.5,-2.4) {};
    \node[fill = none, draw = none] at (1.5,-1.6) {};
    \node[fill = none, draw = none] at (1.5,-.8) {};
\end{tikzpicture}
\end{center}
\caption{The on $1\dot{2}$ pattern stack-sorting map on $\pi=3124$}
\label{fig2}
\end{figure}

In this paper, we introduce another extension of West's stack-sorting map. Namely, we introduce stack-sorting maps $s_{\dot{\tau}}$ that avoid dotted patterns $\dot{\tau}$ (see for example, \Cref{fig2}). Similarly, inspired by Cerbai, Claesson, and Ferrari's \cite{stacksort} $\sigma$-machines, we establish $\dot{\tau}$-machines, which consist of the dotted pattern-avoiding map $s_{\dot{\tau}}$ followed by $s$, and are to be denoted by $s \circ s_{\dot{\tau}}$. We extend Baril's \cite{baril} concept of dotted pattern avoidance to stack-sorting maps, analogous to Goh's \cite{goh2020variations} thesis which introduced barred pattern avoidance to stack-sorting maps. 

Dotted patterns were first introduced by Baril \cite{baril} in 2011. To define dotted patterns, it is first necessary to define weak avoidance of barred patterns. A barred pattern $\bar{\tau}$ is a permutation with some \textit{barred} entries denoted by a bar over the entry. For example, $12\bar{3}45\bar{6}$ is a barred pattern. A permutation \textit{weakly avoids} a barred pattern $\bar{\tau}$ if and only if each occurrence of the pattern $\hat{\tau}$, the unbarred elements of $\bar{\tau}$, in $\pi$ can be extended into an occurrence of the pattern $\tau$ in $\pi$. For example, let $\bar{\tau}=12\bar{3}$. Then, the permutation 1243 weakly avoids $\bar{\tau}$, but 3412 does not.

A dotted pattern $\dot{\tau}$ of length $k$ is defined as a permutation $\tau\in S_n$ with dots over a subset of its entries. A permutation avoids the dotted pattern $\dot{\tau}$ if and only if it weakly avoids the barred pattern $\bar{\tau}^j$ for all $j$ such that the $j$th index does not contain a dot, where $\bar{\tau}^j$ is the barred pattern with a bar over only the $j$th index. For example, let $\dot{\tau}=23\dot{1}$. Then, we have $\bar{\tau}^1=\bar{2}31$ and $\bar{\tau}^2=2\bar{3}1$. It follows that the permutations avoiding $23\dot{1}$ in $S_4$ are 1234, 1342, 2314, and 2341.

In this paper, we prove several analogs of the classical results of West's stack-sorting map for our stack-sorting map that avoids dotted patterns. In \Cref{prelim}, we establish preliminaries. In \Cref{tstacksortable}, we enumerate the $t$-stack-sortable permutations under both dotted pattern-avoiding stacks for $1\dot{2}$ (see \Cref{Th_3-4}) and $2\dot{1}$-avoiding (see \Cref{Th_3-6}) stack for all integers $t \ge 1$. In \Cref{machines}, we characterize and enumerate the permutations that are sorted by the $2\dot{1}$-machine (see \Cref{Th_4-2}) and the fixed points of the $2\dot{1}$-machine (see \Cref{Th_4-4}). In \Cref{highmin}, we show what the maximum number of sorts required is to sort each permutation $\pi\in S_n$ with the $1\dot{2}$-machine (see \Cref{Th_5-2}), and the image of $\pi$ when sorted a number of times near the maximum required sorts (see \Cref{Th_5-4}). In \Cref{futuredirections}, we suggest natural future directions.

\section{Preliminaries}
\label{prelim}

Let $S_n$ be the symmetric group of all permutations of the set $\{1, 2, ..., n\}$. Let $\pi = \pi_1\pi_2 \cdots \pi_n$ be a permutation of length $n$; we denote the length of $\pi$ by $|\pi|$. We also define the \textit{reverse} of $\pi$, $\text{rev}(\pi)$, to be the permutation formed by reversing the order of the entries of $\pi$, or $\pi_n\pi_{n - 1} \cdots \pi_1$. Next, we define $\text{inc}(\pi)$ to be the permutation of length $n$ formed by incrementing each element of $\pi$ by $1$, and the $\text{ins}_i(\pi)$ to be the permutation of length $n+1$ formed by inserting a 1 before the $i$th index of $\text{inc}(\pi)$. For instance, if $\pi=21453$, then $\text{inc}(\pi)=32564$ and $\text{ins}_3(\pi)=321564$.

A \textit{peak} of $\pi$ is an entry $\pi_i$ of $\pi$ such that $\pi_i$ is greater than each of $\pi_1, \pi_2, \ldots, \pi_{i - 1}$. Similarly, a \textit{valley} of $\pi$ is an entry $\pi_i$ of $\pi$ such that $\pi_i$ is less than each $\pi_1, \pi_2, \ldots, \pi_{i - 1}$. For example, if $\pi = 24315$, then the peaks of $\pi$ are $\pi_1 = 2$, $\pi_2 = 4$, and $\pi_5 = 5$, while the valleys of $\pi$ are $\pi_1 = 2$, and $\pi_4 = 1$.

A \textit{peak run} is defined as a maximal sequence of consecutive entries such that the first entry is a peak and no other entry is a peak. Similarly, a \textit{valley run} is defined as a maximal sequence of consecutive entries such that the first entry is a valley and no other entry is a valley. For instance, if $\pi = 24315$, then the peak runs of $\pi$ are $2$, $431$, and $5$, whereas the valley runs of $\pi$ are $243$ and $15$. For a permutation $\pi$, let $\pi = P_1 P_2 \ldots P_{k}$ be the unique partition of $\pi$ into its peak runs. Likewise, let $\pi = V_1 V_2 \ldots V_k$ be the unique partition of $\pi$ into its valley runs.

\subsection{Highly and Minimally Sorted Permutations}

For a set $P \in S_n$ of permutations and any stack-sorting map or machine $s_M$, define the \textit{order} of $P$ with respect to $s_M$, denoted by $\text{ord}_{s_M}(P)$, as the smallest positive integer $k$ such that for any permutation $\pi \in P$, the permutation $(s_M)^k(\pi)$ is a periodic point of $s_M$. 

In 2021, Defant \cite{defant} defined a permutation $S_n$ to be \textit{highly sorted} if it is $t-$sorted for some $t$ close to $n$. In 2023, Choi and Choi \cite{choi} extended Defant's definition by defining the highly sorted permutations of any classical pattern-avoiding map $s_\sigma$ to be the permutations which are $t-$ sorted for $t$ close to $\text{ord}_{s_{\sigma}}(S_n)$. In this paper, we generalize Choi and Choi's definition \cite{choi} by defining the highly sorted permutations of $s_M$ to be the permutations which are $t-$sorted by $s_M$ for $t$ close to $\text{ord}_{s_M}(S_n)$.

As Zhang did in \cite{zhang}, we define a permutation $\pi \in S_n$ to be minimally sorted with respect to $s_M$ if $\text{ord}_{s_M}(\{\pi\}) = \text{ord}_{s_M}(S_n)$. 

For the rest of this paper, we focus on stacks that avoid dotted patterns of length $2$. Of the length-2 dotted patterns, there are 4 possible patterns to study: $\dot{1}2$, $1\dot{2}$, $\dot{2}1$, and $2\dot{1}$. We first show that for dotted patterns of length 2, the placement of the dot is inconsequential.

\begin{theorem} \label{reduction}
    It holds that $s_{1\dot{2}} = s_{\dot{1}2}$ and $s_{\dot{2}1} = s_{2 \dot{1}}$. 
\end{theorem}
\begin{proof} Under $s_{1\dot{2}}$ and $s_{\dot{1}2}$, an entry is pushed onto the stack if and only if it can be extended into an occurrence of the $12$ pattern, so the two stack-sorting maps are identical. Likewise, we have that $s_{\dot{2}1} = s_{2 \dot{1}}$. 
\end{proof}

As result of \Cref{reduction}, for the rest of the paper, we focus on $s_{1 \dot{2}}$ and $s_{2 \dot{1}}$.


\section{Enumeration of $t$-stack-sortable permutations}
\label{tstacksortable}
In this section, we determine the the $t$-stack-sortability of permutations in $S_n$ when passed through the $1\dot{2}$-avoiding stack-sorting map (see \Cref{Th_3-4}) and the $\dot{2}1$-avoiding stack-sorting map (see \Cref{Th_3-6}). We first derive an expression for $s_{1\dot{2}}(\pi)$ in terms of the peak runs of $\pi$.

\begin{proposition}\label{Pr_3-1}
    Let $\pi=P_1P_2\ldots P_k \in S_n$. Then, $s_{1\dot{2}}(\pi)=\text{rev}(P_1)\text{rev}(P_2)\ldots\text{rev}(P_k)$.
\end{proposition}
\begin{proof}
    When any given peak run $P_i$ is first passed through the $1\dot{2}$ stack-sorting map, the following must occur, in order.
    \begin{enumerate}[label=\arabic*.]
        \item Any elements in the stack from previous peak runs are popped, as they must be less than the peak that is at the first index of $P_i$.
        \item Every element of $P_i$ is pushed, because the peak of $P_i$ must be pushed into the empty stack, and every subsequent entry in $P_i$ is less than the peak of $P_i$.
        \item Every element of $P_i$, which is now in the stack, is popped, due to being less than the peak of $P_{i+1}$. If there is no $P_{i+1}$, each element is still popped since no push operations can be completed.
    \end{enumerate}
    Consequently, each peak run is reversed in the output from its order in the input. It follows that $s_{1\dot{2}}(\pi)=\text{rev}(P_1)\text{rev}(P_2)\ldots\text{rev}(P_k)$.
\end{proof}

We next use \Cref{Pr_3-1} to prove that every permutation in $S_n$ is $(n-1)$-stack-sortable.

\begin{lemma} \label{Le_3-2}
    Every permutation in $S_n$ is $(n-1)$-stack-sortable under the $1\dot{2}$ stack-sorting map.
\end{lemma}

\begin{proof}
    We claim that if $\pi_{i}=i$ for some index $i$ and $\pi_{j}=j$ for all $j>i$, the entry $j$ remains in the $j$th position for every subsequent sort for all such $j$. Because each element with an index greater than or equal to $i$ is in its own peak run, by \Cref{Pr_3-1}, each one stays in the same index. Since every element with an index less than $i$ will be popped before $i$ is pushed into the stack, the remaining $i-1$ entries form a permutation of length $i-1$ which is sorted before all other entries are sorted in place. As such, after one sort, the element $\pi_{i-1}$ must equal $i-1$. Each sort causes one element to be added to the ending sequence of sorted entries, and therefore every permutation in $S_n$ is $(n-1)$-stack-sortable under the $1\dot{2}$ stack-sorting map since a permutation of length $n$ is sorted if $n-1$ elements are.
\end{proof}

Now, to invoke induction in our enumeration of $t$-stack-sortable permutations under $s_{1\dot{2}}$, we introduce the following lemma.

\begin{lemma} \label{Le_3-3} 
    The relative order of each element besides 1 in $s_{1\dot{2}}(\text{ins}_i(\pi))$ and each element in $s_{1\dot{2}}(\pi)$ is the same (with an element $x$ not equal to 1 in $(\text{ins}_i(\pi))$ corresponding with $x-1$ in $\pi$). 
\end{lemma}

\begin{proof}
    If $i=1$, then the 1 must be the only entry in the first peak run. By \Cref{Pr_3-1}, $s_{1\dot{2}}(\text{ins}_1(\pi))$ is a 1 followed by $s_{1\dot{2}}(\text{inc}(\pi))$. If 1 is not in the first index, the 1 will be part of some other peak run. It will not create any new peak runs as 1 cannot be a peak if it is not in the 1st index, nor affect the relative order of the other elements in its peak run. By \Cref{Pr_3-1}, the relative order of two elements under $s_{1\dot{2}}$ does not change if in different peak runs and is reversed if in the same peak run. Since inserting a 1 will not affect which peak run any other element is in, the elements 2, 3, \ldots , $n+1$ in $s_{1\dot{2}}(\text{ins}_i(\pi))$ are in the same relative order when compared with the elements 1, 2, \ldots, $n$, in $s_{1\dot{2}}(\pi)$, respectively. Therefore, the relative order of each element besides 1 in $s_{1\dot{2}}(\text{ins}_i(\pi))$ is the same as the relative order of each element in $s_{1\dot{2}}(\pi)$.
\end{proof}

Now, we use the above lemmas to enumerate the $t-$stack-sortable permutations under $s_{1\dot{2}}$.

\begin{theorem} \label{Th_3-4}
    The $t$-stack-sortable permutations under the $1\dot{2}$ stack-sorting map in $S_n$ are enumerated by $n!$ for $n \le t$ and $t! \cdot (t+1)^{n-t}$ otherwise. 
\end{theorem}

\begin{proof}
    If $n\le t$, every permutation in $S_n$ is $(n-1)$-stack-sortable under the $1\dot{2}$ stack-sorting map by \Cref{Le_3-2}. Therefore, the $t$-stack-sortable permutations are enumerated by $n!$ for $n \le t$. If $n > t$, we will prove the $t$-stack-sortable permutations are enumerated by $t! \cdot (t+1)^{n-t}$ by induction on $n$. 
    
    The base case where $n = t + 1$ holds. By \Cref{Le_3-2}, all $n!$ permutations in $S_n$ are $t$-stack-sortable under the $1\dot{2}$ pattern since $t = n - 1$. Indeed, when $n = t + 1$, $t!\cdot(t+1)^{n-t}=t!\cdot(t+1)^{t+1-t}=(t+1)!=n!$. It remains to prove the inductive step. 
    
    By \Cref{Le_3-3}, $s_{1\dot{2}}^t(\text{ins}_i(\pi))$ is equal to the permutation $s_{1\dot{2}}^t\text{inc}(\pi)$ if we remove the 1 from consideration. Therefore, a permutation $\text{ins}_i(\pi)$ is $t$-stack-sortable if and only if the first index of  $s_{1\dot{2}}^t(\pi)$ is 1 and $\pi$ is $t$-stack-sortable. We claim that, given some $\pi\in S_{n-1}$, for the set of $n$ permutations that is all $n$ possible $\text{ins}_i(\pi)$, there are $t+1$ permutations such that $s_{1\dot{2}}^t(\pi)$ is 1. 

    First, observe that the first index of $\text{ins}_1(\pi)$ is 1 after 0 sorts. We claim that after the first sort, the element 1 of each $\text{ins}_i(\pi)$ for $i\ne 1$ gets sorted to a distinct index that is not the $n$th index. Under the $1\dot{2}$ stack-sorting map, the largest element always gets sorted to the last index, so 1 cannot be sorted to the $n$th index. Furthermore, we know by \Cref{Pr_3-1} that peak runs get reversed. Since $i\ne1$, the 1 must be a non-peak element of a peak run, so it is sorted to index $j$, with $j+1$ being the index that $\pi_{i-1}$ is sorted to. Since $\pi_{i-1}$ is distinct for each distinct value of $i$, it follows that $j$, the index of the 1 after the sort, is distinct for each distinct value of $i\ne 1$.
    
    Therefore, the set of $n$ permutations of the form $s_{1\dot{2}}(\text{ins}_i(\pi))$, contains two permutations with a 1 in the 1st index, and one permutation with a 1 in the $k$th index for each $2\le k\le n-1$.  
    
    Note that each permutation in $s_{1\dot{2}}(\text{ins}_i(\pi))$ ends with $n$. Ignoring the $n$, which is already sorted and does not change position with future sorts, the remaining permutations are similar to the set of all $\text{ins}_i(\sigma)$ for some $\sigma\in S_{n-2}$, with the exception that a duplicate exists, as 2 permutations contain a 1 in the first index. We can characterize these permutations ignoring the $n$ as all $\text{ins}_i({\sigma})$ plus an extra $\text{ins}_1({\sigma})$, given that $\sigma$ is equal to $s_{1\dot{2}}(\pi)$ with the element $n$ removed. 
    
    In the same manner, when we sort each $\text{ins}_i({\sigma})$, ignoring $n-1$ and $n$, to get each $\text{ins}_i({\tau})$ with two extra  $\text{ins}_1{\tau}$, for $\tau\in S_{n-3}$ and $\tau$ being $s_{1\dot{2}}^2(\pi)$ without the $n$ and $n-1$. Therefore the number of permutations in $\text{ins}_i^t(\pi)$ with the 1 sorted to the first index is $t+1$.
    
    Since adding a 1 does not affect the relative order of other entries when sorting by \Cref{Le_3-3}, the permutation $\text{ins}_i(\pi)$ is only $t-$stack-sortable if $\pi$ is $t$-stack-sortable and the number of sorts it takes for the 1 to get from the $i$th index to the first index is $t$ or less. Since each permutation in $S_n$ belongs to exactly one set of the form $\text{ins}_i({\pi})$, the number of $t-$stack-sortable permutations of length $n$ is the number of $t-$stack-sortable permutations of length $n-1$ (sets where the incremented elements of $\pi$ get sorted) times $t+1$ (elements in each set where the 1 is sorted). Therefore, incrementing $n$ leads to $t+1$ times as many $t$-stack-sortable permutations. Thus concludes the inductive step for $n$.
\end{proof}

Next, we establish the analog of \Cref{Pr_3-1} for $s_{2\dot{1}}$. 

\begin{proposition}\label{Pr_3-5}
    Let $\pi=V_1V_2\ldots V_k \in S_n$. Then, $s_{2\dot{1}}(\pi)=\text{rev}(V_1)\text{rev}(V_2)\ldots\text{rev}(V_k)$.
\end{proposition}
\begin{proof}
    When any given peak run $V_i$ first passes through the $\dot{2}1$ stack-sorting map, the following must occur, in order.
    \begin{enumerate}[label=\arabic*.]
        \item Any elements in the stack from previous peak runs are popped, as they must be greater than the valley that is at the first index of $V_i$.
        \item Every element of $V_i$ is pushed, because the valley of $V_i$ must be pushed into the empty stack, and every subsequent entry in $V_i$ is greater than the valley of $V_i$.
        \item Every element of $V_i$, which is now in the stack, is popped, due to being greater than the peak of $V_{i+1}$. If there is no $V_{i+1}$, each element is still popped since no push operations can be completed.
    \end{enumerate}
    Consequently, each valley run is reversed in the output from its order in the input. It follows that $s_{\dot{2}1}(\pi)=\text{rev}(V_1)\text{rev}(V_2)\ldots\text{rev}(V_k)$.
\end{proof}

We now use \Cref{Pr_3-5} to enumerate the $t-$stack-sortable permutations for $s_{2\dot{1}}$. 

\begin{theorem} \label{Th_3-6} 
    The $t$-stack-sortable permutations under the $2\dot{1}$ stack-sorting map in $S_n$ are enumerated by $1$ for $n=1$, and 0 otherwise.
\end{theorem}

\begin{proof}
    By \Cref{Pr_3-5}, each valley run is reversed under the $2\dot{1}$ pattern. Since $1$ must be a valley, its index cannot decrease with each sort. Thus, $\pi_{1}=1$, and the valley run containing 1 must be of length 1. It follows that the only $1$-stack-sortable permutation is 1 in $S_1$. It also holds for $t>1$ that the sole $t$-stack-sortable permutation under the $2\dot{1}$ stack-sorting map is $1$.
\end{proof} 

\section{Machine-sortable Permutations and Fixed Points of the $2\dot{1}$-Machine}
\label{machines}
In this section, we characterize and enumerate both the permutations in $S_n$ sorted by the $1\dot{2}$-machine and the fixed points of the $1\dot{2}$-machine. In the following lemma, we characterize the permutations that are sorted by the $\dot{2}1$-machine.  

\begin{lemma} \label{Le_4-1} A permutation $\pi \in S_n$ is sorted under the $2\dot{1}$-machine map if and only if $s_{2\dot{1}}(\pi) = n(n - 1) \cdots 1$. 
\end{lemma}

\begin{proof} In \cite{knuth}, Knuth proved that a permutation is sorted by $s$ if and only if it avoids the pattern $231$ classically. Because $1$ must be a valley of $\pi$, and the start of the last valley run, then the permutation $s_{2\dot{1}}(\pi)$ will end in $1$ by Lemma \ref{Pr_3-5}. If $s_{2\dot{1}}(\pi)$ contains two indices $i, j$ with $1 \le i < j < n$ and $s_{2\dot{1}}(\pi)_i > s_{2\dot{1}}(\pi)_j$, the entries $s_{2\dot{1}}(\pi)_i, s_{2\dot{1}}(\pi)_j$, and $1$ form the classical $231$ pattern in $s_{2\dot{1}}(\pi)$. Therefore $s_{2\dot{1}}(\pi)$ avoids $231$ classically if and only if all of its entries are in decreasing order, or $s_{2\dot{1}}(\pi) = n(n - 1) \cdots 1$.
\end{proof}

Using \Cref{Le_4-1}, we enumerate permutations that are sorted with the $2\dot{1}$-machine.

\begin{theorem} \label{Th_4-2} The machine-sortable permutations under the $2\dot{1}$-machine map in $S_n$ are enumerated by $2^{n - 1}$.
\end{theorem}

\begin{proof} By Lemma \ref{Le_4-1}, we know that $\pi \in S_n$ satisfies $s_{2\dot{1}}(\pi) = n(n - 1) \cdots 1$ exactly when every valley run from an index $i$ to an index $j$ with $i \le j$ satisfies 
$$s_{2\dot{1}}(\pi)_is_{2\dot{1}}(\pi)_{i + 1} \cdots s_{2\dot{1}}(\pi)_j = (n + 1 - j)(n + 2 - j) \cdots (n + 1 - i).$$ Thus, there exists a bijection between the number of ways to split up the $n$ positions of a permutation into valley runs, and length-$n$ permutations sortable under the $2\dot{1}$-machine. Because the first entry must be a valley and all $n - 1$ entries after that must be a valley or in the same valley run as the entry before it, these two sets are both enumerated by $2^{n - 1}$.\end{proof}

 Next, we characterize the fixed points of the $2\dot{1}$-machine. 

\begin{lemma} \label{Le_4-3} Let $\pi =V_1V_2 \cdots V_k \in S_n$. Then all $\pi$ which are fixed points of the $2\dot{1}$-machine are those such that:
    \begin{itemize}
        \item all of the $V_i$ have their entries in increasing order
        \item for $1 \le i < k$, the last entry in $V_{i + 1}$ is larger than all entries of $V_i$
    \end{itemize}   
\end{lemma}

\begin{proof} By Lemma \ref{Pr_3-5}, $s_{2\dot{1}}(\pi)$ must begin with the last entry of $V_1$, or $\pi_{|V_1|}$. In order for $\pi_{|V_1|}$ to be the $|V_1|$th entry popped when $s_{2\dot{1}}(\pi)$  passes through $s$, it must be larger than all of the other entries of $V_1$, and smaller than the last entry of $V_2$ (if $k \ge 1$). Similarly, the second-to-last entry of $V_1$ must be larger than all entries of $V_1$ other than $\pi_{|V_1|}$. Analogous reasoning for the other entries of $V_1$ shows that the entries of $V_1$ must be in increasing order. 

Analogously, we find that all of $V_1, V_2, \ldots, V_k$ have entries in increasing order, and that for $1 \le i < k$, all entries of $V_i$ are lesser than the last entry of $V_{i + 1}$.
\end{proof}

We end this section by using \Cref{Le_4-3} to enumerate the periodic points of the $2\dot{1}$-machine map.

\begin{theorem} \label{Th_4-4} The fixed points of the $2\dot{1}$-machine in $S_n$ are enumerated by OEIS A007476 \cite{oeis}.
\end{theorem}

\begin{proof} We will prove \Cref{Th_4-4} theorem by strong induction on $n$. The base cases of $n = 1$ and $n = 2$ both hold, as $S_1$ and $S_2$ each have one fixed point of the $2\dot{1}$-machine. 

For the inductive step, suppose that $n \ge 3$ and assume that the theorem holds for all of $1, 2, \ldots, n - 1$. Note that the last entry of any fixed point $\pi$ must be $n$, as it will always be sorted to the $n$th position by $s$. Now, we consider cases based on where $n - 1$ lies in $\pi$. 

If $n - 1$ lies in the $n - 1$th position, the last valley run of $\pi$ will end in $(n - 1)n$. By Lemma \ref{Le_4-3}, the first $n - 1$ entries of $\pi$ must be a fixed point of $S_{n - 1}$. There are $a_{n - 1}$ such fixed points, or $\sum_{k = 0}^{n - 3}\tbinom{n - 3}{k}a_k$, where $a_n$ is the $n^{\text{th}}$ entry of OEIS A007476 \cite{oeis}. 

If $n - 1$ lies in the $k$th position for $1 \le k \le n - 2$, it must be immediately followed by a $1$, as the valley run after the one containing $n - 1$ must end in $n$. So, we know that $n$, as the last entry of $\pi$, must be in the valley run beginning with $1$. Then, excluding $1$, $n - 1$, and $n$, suppose we choose $k - 1$ of the remaining $n - 3$ elements to precede $n - 1$; together with $n - 1$, by Lemma \ref{Le_4-3}, we find that the first $k$ entries of $\pi$ form a permutation order-isomorphic to a fixed point of $S_k$. Additionally, the remaining $n - k - 2$ elements between $1$ and $n$ must be in ascending order. Thus, there exists $\binom{n - 3}{k - 1}a_k$ fixed points of $S_n$ with $n - 1$ in the $k$th position. 

Summing over all possible positions of $n - 1$, we find there are
\begin{align*}
    \sum_{k = 0}^{n - 3}\binom{n - 3}{k}a_k + \sum_{k = 1}^{n - 2}\binom{n - 3}{k - 1}a_k &= a_0 + a_{n - 2} + \sum_{k = 1}^{n - 3}\left(\binom{n - 3}{k} + \binom{n - 3}{k - 1}\right)a_k \\ 
    &= a_0 + a_{n - 2} + \sum_{k = 1}^{n - 3}\binom{n - 2}{k}a_k \\
    &= \sum_{k = 0}^{n - 2}\binom{n-2}{k}a_k 
\end{align*}
fixed points of the $2\dot{1}$-machine in $S_n$. The theorem follows from OEIS A007476 \cite{oeis}. 
\end{proof}

\section{Highly and minimally sorted permutations}
\label{highmin}


We begin this section by using \Cref{tstacksortable} to enumerate the highly and minimally sorted permutations of $s_{1\dot{2}}$. As a result of \Cref{Le_3-2} and \Cref{Th_3-4}, we find $\text{ord}_{s_{1\dot{2}}}(S_n) = n - 1$

\begin{corollary} \label{Co_5-1} Under the $1\dot{2}$ stack-sorting map, the number of minimally sorted permutations of length $n$ is $(n-1)!$, and the number of highly sorted permutations is $(n-1)((n-1)!)$. 
\end{corollary}

\begin{proof} The number of $n-1-$stack-sortable permutations under the $s_{1\dot{2}}$ is $n!$ by \Cref{Le_3-2}, and the number of $n-2-$stack-sortable permutations is $(n-2)!\cdot(n-1)^{n-(n-2)}=(n-1)((n-1)!)$ by \Cref{Th_3-4}. Because the identity permutation is a fixed point of the $1\dot{2}$-avoiding stack, it follows that $\text{ord}_{s_{1\dot{2}}}(S_n) = n - 1$. So, the number of highly sorted permutations is $(n - 1)((n - 1)!)$, and the number of minimally sorted permutations is $n!-(n-1)((n-1)!)=(n-1)!$.
\end{proof}

As a follow up to \Cref{Co_5-1}, we characterize image of $s^{n-2}_{1\dot{2}}$.

\begin{theorem} \label{Th_5-2}
    It holds that $s^{n-2}_{1\dot{2}}(S_n) = \{ 1234\ldots (n-1)n, 2134\ldots (n-1)n\}.$
\end{theorem}

\begin{proof} For a permutation $\pi \in S_n$, after $n-2$ sorts under the $1\dot{2}$ pattern, all entries except $1$ and $2$ are necessarily sorted to the correct indices. If 1 is correctly sorted, the permutation is the identity, or $1234\ldots (n-1)n$. Otherwise, the image is $2134\ldots (n-1)n$. Note that both of the images exist for each $n$ since there are nonzero highly and minimally sorted permutations by \Cref{Co_5-1}.
\end{proof}

Next, we find the order of the $1\dot{2}$-machine.

\begin{lemma}  \label{Le_5-3}
For $n \ge 2$, all permutations $\pi \in S_n$ are $\lfloor \tfrac{n}{2} \rfloor $-stack-sortable under the $1\dot{2}$-machine.
\end{lemma}

\begin{proof} We prove \Cref{Le_5-3} by strong induction on $n$. We use $n = 2$ and $n = 3$ as our base cases. The base case of $n = 2$ holds as $(s \circ s_{1\dot{2}})(12) = (s \circ s_{1\dot{2}})(21) = 12$. To verify that the base case of $n = 3$ holds, we note that because $3$ must be the last peak of any permutation $\pi \in S_3$, it must hold that $s_{1\dot{2}}(\pi)$ ends in $3$ by \Cref{Pr_3-1}. Then, when $s_{1\dot{2}}(\pi)$ is passed through $s$, we find that $1$ must be popped from the stack before $2$ enters. Because $(s \circ s_{1\dot{2}})(\pi)$ ends with $3$, we must have $(s \circ s_{1\dot{2}})(\pi) = 123$. 

For the inductive step, assume that $(s \circ s_{1\dot{2}})^{\lfloor k/2 \rfloor}(S_{k}) = 12 \cdots k$ for all $k < n$. For any permutation $\pi \in S_n$, the last peak of $\pi$ must be $n$, meaning $s_{1\dot{2}}(\pi)$ ends in $n$ by \Cref{Pr_3-1}. When $s_{1\dot{2}}(\pi)$ passes through $s$, the elements $1, 2, \ldots, n - 2$ must all be popped from the stack before $n - 1$ enters. Since the last element of $(s \circ s_{1\dot{2}})(\pi)$ is $n$, we find that $(s \circ s_{1\dot{2}})(\pi) = \pi'(n - 1)n$ for some $\pi' \in S_{n - 2}$. When the permutation $(s \circ s_{1\dot{2}})(\pi) = \pi'(n - 1)n$ is repeatedly passed through the $1\dot{2}$-machine, the elements $n - 1$ and $n$ do not change positions. Because $\pi'$ sorts after $\lfloor \tfrac{n - 2}{2} \rfloor$ passes through the $1\dot{2}$-machine by the inductive hypothesis, we see that $\pi$ gets sorted after $\lfloor \tfrac{n - 2}{2} \rfloor + 1 = \lfloor \tfrac{n}{2} \rfloor$ passes through the $1\dot{2}$-machine.
\end{proof}

Naturally, we follow up \Cref{Le_5-3} by characterizing the image of $(s \circ s_{1\dot{2}})^{\lfloor n/2 \rfloor - 1}$. 

\begin{theorem}  \label{Th_5-4} For $n \ge 4$, it holds that
$$(s \circ s_{1\dot{2}})^{\lfloor n/2 \rfloor - 1}(S_n) = \{1234 \cdots (n - 1)n, 2134\cdots (n - 1)n\}$$
if $n$ is even, and 
\begin{align*}
    (s \circ s_{1\dot{2}})^{\lfloor n/2\rfloor - 1}(S_n) = \{&12345 \cdots (n - 1)n, 13245 \cdots (n - 1)n, \\
    &21345 \cdots (n - 1)n, 23145 \cdots (n - 1)n, \\
    &31245 \cdots (n - 1)n\}
\end{align*}
 if $n$ is odd.
\end{theorem}

\begin{proof} First, we prove \Cref{Th_5-4} for even $n$. For such $n$, after $\pi$ goes through the $1\dot{2}$-machine $\lfloor \tfrac{n}{2} \rfloor - 1 = \tfrac{n}{2} - 1$ times, the last $n - 2$ digits must be sorted, so the only possibilities for the resulting permutation are $1234 \cdots (n - 1)n$ and $2134 \cdots (n - 1)n$. Because $1234 \cdots (n - 1)n$ gets sorted to itself by the $1\dot{2}$-machine, it must be in the image after $\tfrac{n}{2} - 1$ sorts. 

To prove that $2134 \cdots (n - 1)n$ is also in the image, we will prove by induction on $n$ that we have $(s \circ s_{1\dot{2}})^{(n/2) - 1}(24 \cdots n 13 \cdots (n - 1)) = 2134 \cdots (n - 1)n$ for all even $n \ge 4$. For the base case of $n = 4$, we note that $s_{1\dot{2}}(2413) = 2314$ and $s(2314) = 2134$. For the inductive step, assume that $(s \circ s_{1\dot{2}})^{((n - 2)/2) - 1}(24 \cdots (n - 2)13 \cdots (n - 3)) = 2134 \cdots (n - 3)(n - 2)$. Then, notice that
\begin{align*}
(s \circ s_{1\dot{2}})^{(n/2) - 1}(24 \cdots n 13 \cdots (n - 1)) &= (s \circ s_{1\dot{2}})^{((n - 2)/2) - 1}(s(s_{1 \dot{2}}(24 \cdots n 13 \cdots (n - 1)))) \\
&= (s \circ s_{1\dot{2}})^{((n - 2)/2) - 1}(s(24 \cdots (n - 2)(n - 1)(n - 3) \cdots 1n)) \\
&= (s \circ s_{1\dot{2}})^{((n - 2)/2) - 1}(24 \cdots (n - 2)13 \cdots (n - 3)(n - 1)n) \\
&= \left((s \circ s_{1\dot{2}})^{((n - 2/2) - 1}(24 \cdots (n - 2)13 \cdots (n - 3)) \right) (n - 1)n \\
&= 2134 \cdots (n - 2)(n - 1)n,
\end{align*}
which finishes the inductive step.

For odd $n \ge 5$, after $\lfloor \tfrac{n}{2} \rfloor - 2 = \tfrac{n - 5}{2}$ sorts of a permutation $\pi \in S_n$ through the $1\dot{2}$-machine, the resulting permutation must have its last $n - 5$ entries sorted. So, after one more sort, the permutation must be of the form $\pi'67 \cdots n$, where $\pi' \in (s \circ s_{1 \dot{2}})(S_5)$. We may verify that $(s \circ s_{1 \dot{2}})(S_5) = \{12345, 13245, 21345, 23145, 31245\}$, which ensures that the only possibilities for permutations in $(s \circ s_{1\dot{2}})^{\lfloor n/2\rfloor - 1}(S_n)$ are the ones in the theorem statement. 

It now remains to show all five permutations are in the image. Note that $12345 \cdots (n - 1)n$ must be in the image as it is sorted to itself by the $1\dot{2}$-machine. To show that $23145 \cdots (n - 1)n$ is in the image, we show that we must always have $(s \circ s_{1\dot{2}})^{\lfloor n/2 \rfloor - 1}(234 \cdots (n - 1)n1) = 23145 \cdots (n - 1)n$ for odd $n \ge 5$ by induction on $n$. The base case of $n = 5$ holds as $s_{1\dot{2}}(23451) = 23415$ and $s(23415) = 23145$. For the inductive step, assume that  $(s \circ s_{1\dot{2}})^{\lfloor (n - 2)/2 \rfloor - 1}(234 \cdots (n - 3)(n - 2)1) = 23145 \cdots (n - 3)(n - 2)$. Then, we have
\begin{align*}
    (s \circ s_{1\dot{2}})^{\lfloor n/2 \rfloor - 1}(234 \cdots (n - 1)n1) &= (s \circ s_{1\dot{2}})^{\lfloor (n - 2)/2 \rfloor - 1} (s(234 \cdots (n - 1)1n)) \\
    &= (s \circ s_{1\dot{2}})^{\lfloor (n - 2)/2 \rfloor - 1} (234 \cdots (n - 2)1(n - 1)n) \\
    &= \left((s \circ s_{1\dot{2}})^{\lfloor (n - 2)/2 \rfloor - 1}(234 \cdots (n - 2)1) \right) (n - 1)n \\
    &= 23145 \cdots (n - 1)n,
\end{align*}
which completes the induction. 

To show that $21345 \cdots (n - 1)n, 13245 \cdots (n - 1)n$, and $31245 \cdots (n - 1)n$ are in the desired image, fix $\pi^3$ to be one of $213$, $132$, and $312$. Define the permutations $\tau^3$ and $\sigma^3$ such that if $\pi^3 = 213$, then $\tau^3 = 2$ and $\sigma^3 = 13$; if $\pi^3 = 132$, then $\tau^3 = 13$ and $\sigma^3 = 2$; and if $\pi^3 = 312$, then $\tau^3 = 3$ and $\pi^3 = 12$. Then, for all odd $n \ge 5$, define the permutations $\pi^n, \tau^n$, and $\sigma^n$ recursively so that $\tau^n$ is $\tau^{n - 2}$ followed by the last entry of $\tau^{n - 2}$ increased by $2$, the permutation $\sigma^n$ is $\sigma^{n - 2}$ followed by the last entry of $\sigma^{n - 2}$ increased by $2$, and $\pi^n = \tau^n\sigma^n$.

We will now prove by induction on $n$ that for all odd integers $n \ge 5$, it holds that $(s \circ s_{1\dot{2}})^{\lfloor n/2 \rfloor - 1}(\pi^n) = \pi^345 \cdots (n - 1)n$. The base case of $n = 5$ holds for all three values of $\pi^3$. For the inductive step, assume that $(s \circ s_{1\dot{2}})^{\lfloor (n-2)/2 \rfloor - 1}(\pi^{n - 2}) = \pi^345 \cdots (n - 3)(n - 2)$. If $\pi_3 \in \{132, 312\}$, we find that 
\begin{align*}
    (s \circ s_{1\dot{2}})(\pi^n) &= (s \circ s_{1\dot{2}})(\tau^{n - 2}n\sigma^n) \\
    &= s(\tau^{n - 2}\text{rev}(\sigma^n)n) \\ 
    &= \tau^{n - 2}\sigma^{n - 2}(n - 1)n \\ 
    &= \pi^{n - 2}(n - 1)n.
\end{align*}If $\pi_3 = 213$, then we have 
\begin{align*}
    (s \circ s_{1 \dot {2}})(\pi^n) &= (s \circ s_{1\dot{2}})(\tau^{n - 2}(n - 1)\sigma^{n - 2}n) \\
    &= s(\tau^{n - 2}\text{rev}(\sigma^{n - 2})(n - 1)n) \\
    &= \tau^{n - 2}\sigma^{n - 2}(n - 1)n \\
    &= \pi^{n - 2}(n - 1)n.
\end{align*} So, no matter what $\pi^3$ is, we obtain
\begin{align*}
    (s \circ s_{1\dot{2}})^{\lfloor n/2 \rfloor - 1}(\pi^n) &= (s \circ s_{1\dot{2}})^{\lfloor (n-2)/2 \rfloor - 1} (\pi^{n - 2}(n - 1)n) \\
    &= \left((s \circ s_{1\dot{2}})^{\lfloor (n-2)/2 \rfloor - 1}(\pi^{n - 2}) \right) (n - 1)n \\
    &= \pi^3 45 \cdots (n - 1)n,
\end{align*}
ending the inductive step.
\end{proof}

\section{Future Directions}
\label{futuredirections}
In \Cref{Th_4-2}, we showed that the number of machine-sortable permutations of $S_n$ under the $2\dot{1}$-machine map is $2^{n - 1}$. We leave the analog of \Cref{Th_4-2} for the $1\dot{2}$-machine as a conjecture. 

\begin{conjecture} \label{Co_6-1}
    The machine-sortable permutations under the $1\dot{2}$-machine map in $S_n$ are enumerated by $\tbinom{2n - 2}{n - 1}$.
\end{conjecture}

We also leave the conjecture that there are no cycles of length $\geq 2$ 
in a $2\dot{1}$-machine.

\begin{conjecture} \label{Co_6-2}
    All permutations in $S_n$ for all $n \ge 1$ are eventually mapped to a fixed point of the $2\dot{1}$-machine after a finite number of iterations through the machine.
\end{conjecture}

Finally, we note that West's stack-sorting map is not limited to only permutations and has also been extended to Coexeter groups by Defant \cite{defant2021stack}, as well as words by Defant and Kravitz \cite{defant2018stack}, and partition diagrams by Campbell \cite{campbell2023lift}. We ask whether dotted patterns could also be generalized to Coexeter groups, words, and partition diagrams, and thus extend dotted pattern-avoiding maps to such mathematical objects. 

\section*{Acknowledgements}
The authors thank Yunseo Choi for suggesting the problem and helpful conversations about the organization of the paper.

\section*{Data Availability Statement}
Data availability not applicable.

\section*{Conflict of Interest}
The authors have no relevant financial or non-financial interests to disclose.


\begin{thebibliography}{20}
\normalsize
\bibitem{avisandnewborn}  D. Avis and M. Newborn, {\em On pop-stacks in series}, Util. Math. 19 (1981).
\bibitem{baril} J. L. Baril, {\em Classical sequences revisited with permutations avoiding dotted pattern}, Electron. J. Combin. 18 (2011).
\bibitem{berlow} K. Berlow, {\em Restricted stacks as functions}, Discrete Math. 344 (2021), 112571.
\bibitem{identitypermutation} M. Bóna, {\em Stack words and a bound for 3-stack sortable permutations}, Discrete Appl. Math. 284 (2020).
\bibitem{campbell2023lift} J. M. Campbell, {\em A lift of West’s stack-sorting map to partition diagrams}, Pacific J. Math. 324.2 (2023) pp. 227-248.
\bibitem{stacksort} G. Cerbai and A. Claesson and L. Ferrari, {\em Stack sorting with restricted stacks}, J. Combin. Theory Ser. A. 173 (2020).
\bibitem{choi} Y. Choi and Y. Choi, {\em Highly sorted permutations with respect to a 312-avoiding stack}, Enumer. Comb. Appl. 3.1 (2022)
\bibitem{defant} C. Defant, {\em Highly sorted permutations and bell numbers}, Enumer. Combin. Appl. 1:1 (2021), $\#$S2R6.
\bibitem{defant2021stack} C. Defant, {\em Stack-sorting for Coxeter groups}, Comb. Theory 2.1 (2021).
\bibitem{defant2018stack} C. Defant and N. Kravitz, {\em Stack-sorting for words}, Australas. J. Combin. 77.1 (2018).
\bibitem{colinandkai} C. Defant and K. Zheng, {\em Stack-sorting with consecutive-pattern-avoiding stacks}, Adv. Appl. Math. 128 (2021).
\bibitem{goh2020variations} Y. K. Goh, {\em Variations of Stack Sorting and Pattern Avoidance}, University of Technology Sydney (Australia) (2020).
\bibitem{knuth} D. E. Knuth, {\em The Art of Computer Programming, Volume 1: Fundamental Algorithms}, Addison-Wesley (1973).
\bibitem{oeis} OEIS Foundation Inc., {\em The On-Line Encyclopedia of Integer Sequences}, (2024) URL: \url{https://oeis.org/}.
\bibitem{west} J. West, {\em Permutations with restricted subsequences and stack-sortable permutations}, Ph.D. Thesis (1990).
\bibitem{xia} J. Xia, {\em Deterministic stack-sorting for set partitions}, Enumer. Comb. Appl. 4.3 (2024).
\bibitem{zeil} D. Zeilberger, {\em A proof of Julian West's conjecture that the number of two-stack-sortable permutations of length $n$ is $2(3n)!/((n+1)!(2n+1)!)$}, Discrete Math. 102 (1992), 85--93. 
\bibitem{zhang} O. Zhang, {\em The Order of the (123, 132)-Avoiding Stack Sort}, Enumer. Comb. Appl. 5.3 (2024).













\end{thebibliography}
\end{document}